\documentclass[a4paper,12pt,twoside]{article}
\usepackage{ifpdf}
\ifpdf
\PassOptionsToClass{pdftex}{article}\relax
\fi
\usepackage{a4wide}
\usepackage{enumerate}
\usepackage{amsmath}
\usepackage{amssymb}
\usepackage{wasysym}
\ifpdf\usepackage[backref]{hyperref}\fi
\usepackage{index}
\usepackage{array}
\usepackage{delarray}
\usepackage{theorem}
\usepackage{comment}
\usepackage[latin1]{inputenc}
\usepackage{xspace}
\usepackage{url}

\newtheorem{thm}{Theorem}

\newtheorem{prop}[thm]{Proposition}
\newtheorem{defn}[thm]{Definition}
\newenvironment{pf}[1][]%
        {\par\noindent{\sc Proof\xspace #1}:~\ignorespaces}%
        {\unskip\nobreak\hskip2em plus1fill$\square$\par\medskip}
\newtheorem{cor}[thm]{Corollary}
\newtheorem{lem}[thm]{Lemma}
\theorembodyfont{\upshape}
\newtheorem{rem}[thm]{Remark}

\makeatletter

\def\@Not#1{\index[notations]{#1}}
\def\@@Not#1{\ifmmode\@Not{#1@{"$#1"$}}\else\@Not{#1}\fi}
\def\Notation#1{\@ifnextchar[{\@Not{#1}}{\@@Not{#1}}}

{\catcode`!\active\global\def!{ }}
\def\@idx[#1]#2{\index{#1}{\catcode`!\active\em #2}}
\def\@@idx#1{\ifmmode\@idx[#1@{"$#1"$}]{#1}\else\@idx[#1]{#1}\fi}
\def\idx{\@ifnextchar[{\@idx}{\@@idx}}%]

\def\ifnextchar{\@ifnextchar}

\let\pointdexclamation=!
\let\inferieur=<
\let\superieur=>
\let\egal==
{
\catcode`!\active
\catcode`<\active
\catcode`>\active
\global\everymath={\catcode`!=\active\def!{\mathpointdexclamation}\catcode`<\active\catcode`>\active%
\def\mathpointdexclamation{\@ifnextchar\egal{\not}{\pointdexclamation}}%
\def<{\@ifnextchar\egal{\leqslant\@gobble}{\inferieur}}%
\def>{\@ifnextchar\egal{\geqslant\@gobble}{\superieur}}}
\global\everydisplay={\catcode`!\active\catcode`<\active\catcode`>\active%
\def!{\@ifnextchar\egal{\not}{\pointdexclamation}}%
\def<{\@ifnextchar\egal{\leqslant\@gobble}{\inferieur}}%
\def>{\@ifnextchar\egal{\geqslant\@gobble}{\superieur}}}
}

\def\aquot#1#2{\mathchoice%
{\leavevmode\kern-.1em\lower.25ex\hbox{$\displaystyle #1$}\kern-.1em\backslash\kern-.1em%
\raise.2ex\hbox{$\displaystyle #2$}}%
{\leavevmode\kern-.1em\lower.25ex\hbox{$\textstyle #1$}\kern-.1em\backslash\kern-.1em%
\raise.2ex\hbox{$\textstyle #2$}}%
{\leavevmode\kern-.1em\lower.25ex\hbox{$\scriptstyle #1$}\kern-.1em\backslash\kern-.1em%
\raise.2ex\hbox{$\scriptstyle #2$}}%
{\leavevmode\kern-.1em\lower.25ex\hbox{$\scriptscriptstyle #1$}\kern-.1em\backslash\kern-.1em%
\raise.2ex\hbox{$\scriptscriptstyle #2$}}}
\def\quot#1#2{\mathchoice%
{\leavevmode\kern-.1em\raise.2ex\hbox{$\displaystyle #1$}\kern-.1em/\kern-.1em%
\lower.25ex\hbox{$\displaystyle #2$}}%
{\leavevmode\kern-.1em\raise.2ex\hbox{$\textstyle #1$}\kern-.1em/\kern-.1em%
\lower.25ex\hbox{$\textstyle #2$}}%
{\leavevmode\kern-.1em\raise.2ex\hbox{$\scriptstyle #1$}\kern-.1em/\kern-.1em%
\lower.25ex\hbox{$\scriptstyle #2$}}%
{\leavevmode\kern-.1em\raise.2ex\hbox{$\scriptscriptstyle #1$}\kern-.1em/\kern-.1em%
\lower.25ex\hbox{$\scriptscriptstyle #2$}}}
\def\biquot#1#2#3{\mathchoice%
{\leavevmode\kern-.1em\lower.25ex\hbox{$\displaystyle #1$}\kern-.1em\backslash\kern-.1em%
\raise.2ex\hbox{$\displaystyle #2$}\kern-.1em/\lower.25ex\hbox{$\displaystyle #3$}}%
{\leavevmode\kern-.1em\lower.25ex\hbox{$\textstyle #1$}\kern-.1em\backslash\kern-.1em%
\raise.2ex\hbox{$\textstyle #2$}\kern-.1em/\kern-.1em\lower.25ex\hbox{$\textstyle #3$}}%
{\leavevmode\kern-.1em\lower.25ex\hbox{$\scriptstyle #1$}\kern-.1em\backslash\kern-.1em%
\raise.2ex\hbox{$\scriptstyle #2$}\kern-.1em/\kern-.1em\lower.25ex\hbox{$\scriptstyle #3$}}%
{\leavevmode\kern-.1em\lower.25ex\hbox{$\scriptscriptstyle #1$}\kern-.1em\backslash\kern-.1em%
\raise.2ex\hbox{$\scriptscriptstyle #2$}\kern-.1em/\kern-.1em\lower.25ex\hbox{$\scriptscriptstyle #3$}}}
\makeatother
\def\oper#1{\ensuremath{\mathop{\textup{#1}}}\xspace}
\def\defoper#1{\expandafter\def\csname #1\endcsname{\oper{#1}}}

\defoper{diag}
\defoper{dim}
\def\FM{\textbf{\upshape F}\xspace}
\def\Fp{{\ensuremath{\FM_p}}\xspace}
\defoper{Frob}
\defoper{Gal}
\def\GL{\textup{GL}\xspace}
\defoper{Gram}

\def\Hom{\oper{Hom}\nolimits}

\defoper{Im}

\def\L{\textup{L}\xspace}
\def\MC{{\mathcal M}\xspace}

\def\NM{\textbf{\upshape N}\xspace}

\def\O{{\ensuremath{\mathcal O}}\xspace}

\def\pG{{\ensuremath{\mathfrak p}}\xspace}

\def\QM{\textbf{\upshape Q}\xspace}
\def\Qp{{\ensuremath{\QM_p}}\xspace}
\defoper{Re}

\defoper{Sp}
\defoper{Stab}

\defoper{Tr}

\defoper{vol}
\defoper{Vect}

\def\ZM{\textbf{\upshape Z}\xspace}

\allowdisplaybreaks

\begin{document}
\nocite{PARI2}
\title{Comparison of semi-simplifications\\of Galois representations}
\author{Lo\"\i c Greni\'e\\
Dipartimento di Matematica ``Felice Casorati''\\Università degli Studi di Pavia\\via Ferrata, 1\\I-27100 Pavia\\Italy\\
loic.grenie@gmail.com}
\maketitle

\section*{Introduction}
The aim of this paper is to provide a criterion to ensure that the
semi-simplifications of $p$-adic finite dimensional Galois representations
are isomorphic. Such an isomorphism implies that the Artin \L-functions of
these representations are the same. This can for example be used to compare the
Artin \L-functions obtained from automorphic representations with those issued
from algebraic geometry. Another possible use is shown in Section
\ref{comparison} which proves that the two representations considered in
\cite{GeemTop:non-self-dual}, one of which is a subrepresentation of the
cohomology of a variety while the other is a conjectural automorphic
representation, are isomorphic. The main result of this paper, Theorem
\ref{theo}, provides an effective criterion to check whether two
semi-simplifications are isomorphic and explains which Frobenius elements
suffice to compare the two. In
\cite[section 4]{Livne:cubic-exponential}, Ron Livné explained and generalized
(by lowering the required number of comparisons) a result of Jean-Pierre Serre
giving a sufficient condition for semi-simplifications of $p$-adic Galois
representations to be isomorphic. We intend to generalize here the original
result of Jean-Pierre Serre. Even though our result is valid for all dimensions
and cannot use the fact that a group of exponent $2$ is abelian, our result is
similar in complexity to the one of Livné\footnote{Livné has a better result
because he shows that he does not need to compare the representations for all
Frobenius elements but only for a so called ``non cubic'' family}.

This paper can be considered as an application of the method explained in
\cite[p27--29]{Serre:oeuvres-4}.

I would like to thank Bert van Geemen who brought this subject to my attention,
gave me some hints and helped remove some errors. This result has been made
possible thanks to the help of Thomas Weigel who helped me understand the
complexity of the pro-$p$ groups and suggested the use of the ``powerful
pro-$p$ groups'' which proved to be the right way to tame that complexity. I
also would like to thank Karim Belabas who has been of great help for the
computational part of the result. I would finally like to thank the referee for
his extremely careful reading and his very precise and insightful remarks.
\section{The result}
\subsection{Setup}
This section sets up the framework of this work. Let us fix an integer $n>=2$,
a prime $p$ and define $m$ as the minimum integer such that $p^m>=n$. We fix a
global field $K$ and let $\overline K$ be a maximal separable algebraic
extension of $K$. All extensions of $K$ considered in this paper are
sub-extensions of $\overline K$. For any subfield $L$ of $\overline K$, we
denote $\Gamma_L=\Gal\left(\quot{\overline K}L\right)$.

We denote $\MC(n,A)$ the algebra of matrices of size $n\times n$ with
coefficients in a ring $A$.
\begin{defn}\label{congruent}
Let $E$ be a finite extension of \Qp for some prime $p$. Let $M_1$ and $M_2$ be
two matrices in $\MC(n,E)$. Let $F$ be a finite extension of $E$ containing the
eigenvalues of $M_1$ and $M_2$. Denote $\O_F$ the integer ring of $F$, $\pG_F$
its maximal ideal and $\varpi_F$ a uniformiser. The two matrices $M_1$ and
$M_2$ are said to have congruent eigenvalues if there exist
$\lambda\in\O_F^{\,\times}$ and $v\in\ZM$ such that the characteristic
polynomials of $\varpi_F^{\,-v}M_1$ and $\varpi_F^{\,-v}M_2$ are in $\O_F[X]$
and are congruent to $(X-\lambda)^n$ modulo $\pG_F$.
\end{defn}
\begin{rem}
\begin{itemize}
\item The absolute Galois group of a global field is compact, so the
eigenvalues of the matrices in the image of a Galois representation
necessarily have valuation $v=0$.
\item The condition on the matrices is rather strong: it implies that the $2n$
eigenvalues are congruent to a \underline{single} one. Therefore the condition
is strong even for each matrix separately.
\end{itemize}
\end{rem}
\subsection{Construction}\label{K_S}
For any finite set of places $S$ of $K$, we want to construct an extension
$K_S=K_{S,n}$ of $K$ such that the Galois group of $\quot{K_S}K$ is sufficient
to compare the semi-simplifications of representations of $\Gamma_K$ with
values in $\GL(n,E)$, unramified outside $S$, and with all eigenvalues reducing
to a single one in the residual field of $E$.

Take $K_0=K$.
Define $K_i$ by induction by taking $K_{i+1}$ to be the maximal abelian
extension of $K_i$ unramified outside $S$ and such that
$\Gal(\quot{K_{i+1}}{K_i})$ is a direct product of copies of \Fp.
Notice that $K_i$ is a Galois extension of $K$ at each step $i$.
Let $\epsilon=0$ if $p!=2$ and $\epsilon=1$ if $p=2$. Let
$r=N^{2(1+\epsilon)}\frac{N(N-1)}2$ with $N=n[E:\Qp]$.
Let $\lambda$ be the minimum integer such that $2^\lambda>=r$.
Finally take $K_S=K_{\lambda+\epsilon+m}$.
\subsection{Main result}
\begin{thm}\label{theo}
We fix an integer $n>=2$, a prime $p$ and define $m$ to be the minimum integer
such that $p^m>=n$.  Let $K$ be a global field, $S$ a finite set of places of
$K$ and $E$ a finite extension of \Qp. We assume that if $k$ is the residual
field of $E$, then $n$ and $|k^\times|$ are relatively prime. Let $K_S$ be the
field constructed as in Section \ref{K_S}. Fix a set $T$ of places of $K$,
disjoint from $S$, such that each maximal cyclic subgroup of
$\Gal(\quot{K_S}K)$ has a generator of the form $\Frob(\quot{\mathfrak t}t)$
for some $t\in T$ and some prime $\mathfrak t$ above $t$ in $K_S$. Assume now
that
$$\rho_1,\rho_2:\Gamma_K\longrightarrow\GL(n,E)$$
are continuous representations unramified outside $S$ and satisfy the following
conditions:
\begin{enumerate}
\item\label{all-eigenvalues-equal}$\forall\sigma\in\Gamma_K$,
$\rho_1(\sigma)$ and $\rho_2(\sigma)$ have congruent eigenvalues (see
Definition \ref{congruent}).
\item \label{cond-polchar}
$\forall t\in T$, $\rho_1(\Frob t)$ and $\rho_2(\Frob t)$ have equal
characteristic polynomials (where $\Frob t$ is any Frobenius element above
$t$).
\end{enumerate}
Then $\rho_1$ and $\rho_2$ have isomorphic semi-simplifications.
\end{thm}
\begin{rem}\begin{itemize}
\item
In this theorem, the condition $(n,|k^\times|)=1$ is needed just to ensure
that, up to a twist by a character, the residual representations are
$p$-groups.
\item
The characteristic polynomial of a matrix $M$ of size $n$ has coefficients
which are symmetric functions of degree at most $n$ of the eigenvalues of $M$.
Over a field of characteristic $0$ the sums of the powers of the variables
$(X_i)_{1<=i<=n}$ are a basis of the space of symmetric functions in $(X_i)$.
It follows that there exists a function $f$ independent of $M$ such that the
characteristic polynomial of $M$ is $f(\Tr M,\Tr M^2,\dots,\Tr M^n)$. Hence we
can modify Condition $(\ref{cond-polchar})$ above as follows:

\noindent{}either
$$\forall t\in T,\forall 1<=k<=n,\Tr\rho_1((\Frob t)^k)=\Tr\rho_2((\Frob t)^k)$$
or
$$\forall t\in T,
\begin{cases}
\forall 1<=k<=n-1,	& \Tr\rho_1((\Frob t)^k)=\Tr\rho_2((\Frob t)^k) \\
			& \det\rho_1(\Frob t)=\det\rho_2(\Frob t).
\end{cases}$$
\item As for the condition $(n,|k^\times|)=1$: observe that, if $n$ is even,
then $p$ has to be $2$. Observe also that if $n$ is a power of $p$, or
$k=\FM_2$, then the condition is verified. Finally observe that we can multiply
$n$ by $[E:\Qp]$ and asume $E=\Qp$. We can thus always apply the theorem if we
choose $p=2$ (at the cost of enlarging $n$, which makes it less interesting
because $K_S$ and $T$ become larger).
\end{itemize}
\end{rem}
\section{Pro-$p$-groups}
The main result of this section is Proposition \ref{central-prop} which
establishes our result for a pro-$p$ group.
\subsection{The result for pro-$p$ groups}
\begin{defn}
For a $p$-group or pro-$p$ group $G$, we denote by $G^\#$ the closure of the
intersection of the kernels of all group morphisms from $G$ to finite groups
such that all their elements have order dividing $p^m$.
\end{defn}
\begin{rem}
\begin{enumerate}
\item $G^\#$ is also called $\textrm{\agemO}_m(G)$, at least when $G$ is
finite.
\item $G^\#$ is normal in $G$.
\item Observe that $G^\#$ is also the closure of the subgroup generated by
$p^m$-th powers.
\item In case $n=p=2$, the subgroup $G^\#$ is just the Frattini subgroup $G^*$.
\item If $\rho:G\longrightarrow H$ is a continuous group morphism,
$\rho(G^\#)=\rho(G)^\#\subseteq H^\#$ with equality if $\rho$ is surjective.
\item If $G_1$ and $G_2$ are groups, then
$(G_1\times G_2)^\#=G_1^\#\times G_2^\#$.
\end{enumerate}
\end{rem}
The following lemma will be useful later on:
\begin{lem}\label{Gamma-bar-cancelletto-non-c'e}
Let $G$ be a $p$-group such that any element of $\quot G{G^\#}$ has a
representative in $G$ of order dividing $p^m$. Then $G^\#=\{1\}$.
\end{lem}
\begin{pf}
Suppose that $G^\#!=\{1\}$. Observe first that, according to
\cite[Theorem 1.12, p90]{Suzuki1}, we can find a normal subgroup $N$ of $G$
which is a subgroup of index $p$ of $G^\#$.  Then
$\left(\quot GN\right)^\#\simeq\quot{G^\#}N\simeq\Fp$, so that we can as well
assume that $G^\#=\Fp$. We have an exact sequence
$$0\longrightarrow\Fp\longrightarrow G\longrightarrow\quot G{G^\#}\longrightarrow 1$$
such that each element of $\quot G{G^\#}$ has a representative in $G$ of order
dividing $p^m$.

Denote $H=\quot G{G^\#}$. Then $H$ is a $p$-group and $\oper{Aut}\Fp$ has $p-1$
elements, so that the action of $H$ on \Fp is trivial. This means that
the extension
$$0\longrightarrow\Fp\longrightarrow G\mathrel{\mathop{\longrightarrow}^\pi}H\longrightarrow 1$$
is central, i.e. that $G^\#\subseteq Z(G)$. Thus every element $g$ of $G$ has
order dividing $p^m$ (all elements of $gG^\#$ have the order of $g$, except if
$g\in G^\#$, in which case all elements have order either $p$ or $1$). We
deduct that the identity is a morphism from $G$ to a group having elements of
order dividing $p^m$.  This means that $G^\#=\{1\}$, which is impossible.
\end{pf}
\begin{rem}
This lemma is a generalization of
\cite[Lemma 4.5, p257]{Livne:cubic-exponential}. The definition of $G^\#$
accounts for Remark 4.6.a. below the proof of the lemma in \emph{loc. cit.}.
\end{rem}
\begin{prop}\label{central-prop}
Let $G$ be a pro-$p$ group which is topologically finitely generated and let
$E$ be a finite extension of \Qp. Recall that the integer $m$ used to define
$G^\#$ is the minimum integer such that $p^m>=n$. Assume
$$\rho_1,\rho_2:G\longrightarrow\GL(n,E)$$
are continuous representations and $\Sigma\subset G$ is a subset satisfying:
\begin{enumerate}
\item the image of $\Sigma^\bullet=\{\sigma^k/\sigma\in\Sigma,k\in\NM\}$ in
$\quot G{G^\#}$ is equal to $\quot G{G^\#}$;
\item\label{cond-polchar-prop}
$\forall\sigma\in\Sigma$, $\rho_1(\sigma)$ and $\rho_2(\sigma)$ have the same
characteristic polynomial.
\end{enumerate}
Then $\rho_1$ and $\rho_2$ have isomorphic semi-simplifications.
\end{prop}
\begin{pf}
Let \O be the integer ring of $E$. Since $G$ is compact, it preserves a
full lattice in $E^n$ when acting via each $\rho_i$, for $i=1$, 2. Since
\O is a discrete valuation ring, such a lattice is free over \O. Hence we
may assume $\rho_i(G)\subset\GL(n,\O)$ for each $i=1$, $2$.

Let \pG be the maximal ideal of \O and set $k=\quot\O\pG$. The reduction modulo
\pG of $\rho_i(G)$ is a $p$-group in $\GL(n,k)$. A $p$-Sylow subgroup for
$\GL(n,k)$ is the subgroup of upper triangular unipotent matrices. We can thus
suppose, up to a base change in the lattices above, that the reduction of
$\rho_i(G)$ modulo \pG is included in this subgroup. In particular, for any $g$
in $G$, $(\rho_i(g)-I_n)^n\equiv0\pmod\pG$. We also have that
$\rho_i(g)^{p^m}\equiv I_n\mod\pG$ (in fact we can substitute $p^m$ by any
power of $p$ that is at least equal to the nilpotency order of the reduction
mod \pG of $\rho_i(g)-I_n$).

Now let $M_n=\MC(n,\O)$. We define $\rho:G\longrightarrow M_n\times M_n$ to
be the map $\rho(g)=(\rho_1(g),\rho_2(g))$. Set $M$ to be the linear \O-span of
$\rho(G)$ in $M_n\times M_n$. Then $M$ is an \O-algebra spanned (as an
\O-module) by $\Gamma=\rho(G)$. Let $R=\quot{M}{\pG M}$ and for $g\in G$, we
will denote the image of $\rho(g)$ in $R$ by $\overline g$. Set
$\overline\Gamma=\{\overline g/g\in G\}$. Then $R$ is a $k$-algebra with unity
$\overline1=(I_n,I_n)\bmod{\pG M}$ and spanned by $\overline\Gamma$ as a
$k$-vector space.

We would like to prove that $R$ is spanned over $k$ by
$\overline{\Sigma^\bullet}=\{\overline{\sigma^k}/\sigma\in\Sigma,k\in\NM\}$. We
claim that, for any $\sigma\in\Sigma$, we have
$(\overline\sigma-\overline1)^n=\overline0$ and
${\overline\sigma}^{p^m}=\overline1$. Both these equalities generalize
${\overline\sigma}^2=\overline1$ for $p=n=2$. The point is that equalities in
$\GL(n,k)$ can sometimes be translated to equalities in $R$. Let us first
observe that the characteristic polynomial of $\rho_i(\sigma)\bmod\pG$ is
$(X-1)^n$. This polynomial is the reduction modulo \pG of the characteristic
polynomial of $\rho_i(\sigma)$. Let $\sum_{r=0}^nc_{r,i}X^r$ be the
characteristic polynomial of $\rho_i(\sigma)$. Let
$a_{r,i}=(-1)^{n-r}\binom nr-c_{r,i}$. Then
$(\rho_i(\sigma)-I_n)^n=\sum_{r=0}^na_{r,i}\rho_i(\sigma)^r$ and all
$a_{r,i}\in\pG$. From Hypothesis \ref{cond-polchar-prop}, we know that the
characteristic polynomials are equal and thus $a_{r,1}=a_{r,2}=a_r$. We can
deduct that
{\def\r{\rho_}\let\s=\sigma
\begin{align*}
(\rho(\s)-(I_n,I_n))^n	& = ((\r1(\s)-I_n)^n,(\r2(\s)-I_n)^n)		\\*
		& = (\sum_{r=0}^na_r\r1(\s)^r,\sum_{r=0}^na_r\r2(\s)^r)	\\*
			& = \sum_{r=0}^na_r(\r1(\s)^r,\r2(\s)^r)	\\*
			& = \sum_{r=0}^na_r(\r1(\s),\r2(\s))^r		\\*
			& = \sum_{r=0}^na_r\rho(\s)^r			\\*
			& \in\pG M
\end{align*}}
Thus $N(\sigma)=\overline\sigma-\overline1$ is nilpotent of order (at most)
$n$. This means that for any $r$,
${\overline\sigma}^r=(\overline1+N(\sigma))^r$ is a polynomial in $N(\sigma)$
of degree at most $n-1$. For $r=p^m>=n$, then $\binom{p^m}r$ will be in
$p\ZM\subset\pG$ for all $i\in[1;p^m-1]$. Thus
${\overline\sigma}^{p^m}=\overline1$ (as above we can substitute $p^m$ by any
power of $p$ that is at least equal to the nilpotency order of each
$\rho_i(\sigma)-I_n\mod\pG$). In addition, since we have only used the fact
that $\rho_1(\sigma)$ and $\rho_2(\sigma)$ have the same characteristic
polynomial, this remains true for all powers of all the elements of $\Sigma$:
$$\forall\sigma\in\Sigma^\bullet,\quad
\begin{cases}
\left(\overline{\sigma}-\overline1\right)^n=\overline0	\\
\overline{\sigma}^{p^m}=\overline1,
\end{cases}
$$
which means
$$\forall\overline\sigma\in\overline{\Sigma^\bullet},\quad
\begin{cases}
(\overline\sigma-\overline1)^n=\overline0	\\
\overline\sigma^{p^m}=\overline1.
\end{cases}
$$

\medskip

To prove that $R$ is $k$-spanned by $\overline{\Sigma^\bullet}$ we first prove
that $\overline\Gamma^\#=\{\overline1\}$. Observe that, since \O is a principal
domain, $R$ is a finite-dimensional $k$-vector space of dimension at most
$2n^2$. Hence $R$ and $\overline\Gamma$ are finite. We can apply Lemma
\ref{Gamma-bar-cancelletto-non-c'e} to show that $\overline\Gamma^\#=\{1\}$:
since $\rho(G)=\Gamma$, we have $\Gamma^\#=\rho(G)^\#=\rho(G^\#)$, which implies
${\overline\Gamma}^\#=\overline{\Gamma^\#}=\overline{\rho(G^\#)}=\overline{G^\#}$
and thus any element of $\quot{\overline\Gamma}{\overline\Gamma^\#}$ can be
represented by an element of $\overline{\Sigma^\bullet}$ and these elements
have order dividing $p^m$. According to Lemma
\ref{Gamma-bar-cancelletto-non-c'e}, we have $\overline\Gamma^\#=\{1\}$ and thus
$\overline\Gamma\simeq\quot{\overline\Gamma}{\overline\Gamma^\#}\simeq\overline{\quot\Gamma{\Gamma^\#}}\simeq\overline{\quot G{G^\#}}\subseteq\overline{\Sigma^\bullet}$
(the last inclusion is up to the canonical projection from $G$ to
$\quot G{G^\#}$); since $\overline{\Sigma^\bullet}\subseteq\overline\Gamma$ and
both are finite, we conclude that $\overline{\Sigma^\bullet}=\overline\Gamma$.

\medskip

Using the former argument, we can apply Nakayama's lemma to see that
$\Sigma^\bullet$ generates $M$ as an \O-module. Since the characteristic
polynomials of $\rho_1(\sigma)$ and $\rho_2(\sigma)$ are equal, the traces of
$\rho_1(\sigma^k)$ and $\rho_2(\sigma^k)$ are equal for all $\sigma\in\Sigma$
and all $k\in\NM$.  Thus the linear form $\alpha$ on $M$ defined by
$\alpha(a,b)=\Tr a-\Tr b$ is trivial on a generating set of $M$ and thus on all
of $M$. As a consequence, the characteristic polynomials of $\rho_1(g)$ and
$\rho_2(g)$ are equal for all $g\in G$.
\end{pf}
\subsection{Structure of pro-$p$ groups}\label{pro-p-groupes}
A good reference for the following is \cite{Analytic-pro-p-groups}, and in
particular chapter 3.
\begin{defn}
A powerful pro-$p$ group is a pro-$p$ group $G$ such that $\quot G{G^p}$ (resp.
$\quot G{G^4}$ if $p=2$) is abelian, where $G^p$ (resp. $G^4$) is the subgroup
generated by $p$-th (resp. fourth) powers of elements of $G$.
\end{defn}
\begin{prop}
For each finitely generated pro-$p$ group $G$ with a powerful open subgroup,
there is a number $r$ such that any subgroup of $G$ has at most $r$ generators.
\end{prop}
\begin{defn}
The minimal number $r$ above is called the rank of the pro-$p$ group $G$.
\end{defn}
For any integer $r>=1$ we define the integer $\lambda(r)$ as the minimum
$\ell$ such that $2^\ell>=r$.

A proof of the following result is included in the proof of
\cite[Theorem 3.10]{Analytic-pro-p-groups}.
\begin{thm}
\label{powerful-pro-p}
For any pro-$p$ group $G$ of rank $r$, there exists a $t<=\lambda(r)+\epsilon$
and a filtration
$$G_t\subseteq G_{t-1}\subseteq\dots\subseteq G_0=G$$
with abelian quotients of exponent $p$ such that $G_t$ is powerful. Recall
that $\epsilon=1$ if $p=2$ and $\epsilon=0$ otherwise.
\end{thm}
\section{Reinterpretation of $\quot G{G^\#}$ in the Galois group}
\begin{pf}[of Theorem \ref{theo}]
We take $k$ to be the residual field of $E$ and $q=|k|$. Since $(n,q-1)=1$, the
map $x\mapsto x^n$ is injective and thus surjective and bijective in $k$.
Therefore there exists a unique character
$$\overline\chi:\Gamma_K\longrightarrow k^\times$$
satisfying ${\overline\chi}^n=\det\rho_i\pmod\pG$ for $i=1$, $2$. Let $\chi$ be
the Teichmüller lift of $\overline\chi$. Then all the eigenvalues of
$\chi^{-1}(g)\rho_i(g)$ will be in some finite extension $F$ of $E$ and they
will reduce to the same $\lambda$ in some finite extension $k'$ of $k$. The
characteristic polynomial of the reduction mod \pG of each
$\chi^{-1}(g)\rho_i(g)$ will be of the form $P_i(X)=(X-\lambda)^n$. We write
$n=p^vm$ with $(m,p)=1$. We then have
$$P(X)=(X^{p^v}-\lambda^{p^v})^m=X^n-m\lambda^{p^v}X^{n-p^v}+...+\lambda^n$$
so that, since $m!=0$ in $k$, $\lambda^{p^v}\in k$. This shows that
$\lambda\in k$. Since $\lambda^n=\overline{\det\rho_i(g)\chi^{-1}(g)}=1$, we
obtain $\lambda=1$. Thus the image of $\Gamma_K$ under the map
$\rho(g)=\chi^{-1}(g)(\rho_1(g),\rho_2(g))$ is a pro-$p$ group
$G\subset\GL(n,E)^2$. This can easily be seen from
\cite[Proposition 1.11, p22]{Analytic-pro-p-groups}: change basis so that both
reductions mod \pG of $\rho_i(\Gamma_K)$, for $i=1$, $2$, have image in the
subgroup $U_k$ of unipotent upper triangular matrices. Let $U_\pG$ be the
inverse image of $U_k$ in $\GL(n,\O)$. Then $U_k$ is a $p$-group and the kernel
of the reduction mod \pG is the normal subgroup $V=I_n+\varpi\MC(n,\O)$ which
is a pro-$p$ group. Hence $U_\pG$ is a pro-$p$ group and $\rho_1(\Gamma_K)$ and
$\rho_2(\Gamma_K)$ are closed in $U_\pG$, because they are compact, thus they
also are pro-$p$ groups.

We want to compute the ranks of $G_1=\rho_1(\Gamma_K)$ and
$G_2=\rho_2(\Gamma_K)$. We begin by embedding $\GL(n,\O)$ in $\GL(N,\ZM_p)$ by
using a basis of $\O$ over $\ZM_p$ to identify $\O^n$ and $\ZM_p^{\,N}$. Let
$M$ be a matrix of $\GL(n,\O)$ with characteristic polynomial $P(X)$. Its
embedding $M_r$ in $\GL(N,\ZM_p)$ has characteristic polynomial equal to
$\prod P^\sigma(X)$, where $\sigma$ runs over the embeddings of $E$ in a fixed
algebraic closure of $E$ and $P^\sigma$ is the polynomial obtained from $P$
by applying $\sigma$ to its coefficients. In particular, if $M$ reduces to an
unipotent matrix in $\GL(n,k)$, its characteristic polynomial is congruent to
$(X-1)^n$ modulo \pG so that the characteristic polynomial of $M_r$ is
congruent to $(X-1)^N$ modulo $p$. This means that $M_r$ reduces to an
unipotent matrix in $\GL(N,\FM_p)$. The group of unipotent matrices of
$\GL(N,\FM_p)$ has rank at most $\frac{N(N-1)}2$. The kernel of the reduction
mod $p$ in $\GL(N,\ZM_p)$ is $V=I_n+p\MC(N,\ZM_p)$. According to \cite[Theorem
5.2, p88]{Analytic-pro-p-groups}, if $p$ is odd then $V$ is powerful of rank
$N^2$ while if $p=2$ then the subgroup $V'=I_n+4\MC(N,\ZM_2)$ is powerful of
rank $N^2$ and $\quot V{V'}$ is a subgroup of $(\quot\ZM{2\ZM})^{N^2}$, which
means that it is a $2$-group of rank at most $N^2$. Putting all three terms
together, we see that the group of matrices that
reduce to the subgroup of unipotent matrices in $\GL(n,k)$ has rank at most
$r=N^2\cdot(N^2)^\epsilon\cdot\frac{N(N-1)}2$. This means that the ranks of
$G_1$ and $G_2$ are at most $r$. We can apply Theorem \ref{powerful-pro-p} to
$G_i$, for $i=1$, $2$: for some $t<=\lambda(r)+\epsilon$, we get a filtration
$$V_i=G_{i,t}\subset G_{i,t-1}\subset\cdots\subset G_{i,1}\subset G_{i,0}=G_i$$
with all quotients $\quot{G_{i,s}}{G_{i,s+1}}$ abelian of exponent $p$ and
$V_i$ a powerful pro-$p$ group. Since $V_i$ is powerful, with $m$ more
filtration steps we get $V_i^\#$. It is clear that since $V_i\subseteq G_i$, we
have $V_i^\#\subseteq G_i^\#$. Since $G^\#$ is the closure of the subgroup
generated by the $p^m$-th powers, we see that
$G^\#\subseteq G_1^\#\times G_2^\#$. This means that a filtration with at most
$\lambda(r)+\epsilon+m$ steps is sufficient to get a subgroup $V^\#$ of
$G^\#$\footnote{Observe that $r$ is not an upper bound for the rank of $G$: the
rank of $G$ is at most $2r$ but can be greater than $r$.}.

On the field side, the $i$-th step of the filtration corresponds to an
extension of $K_i$ by an abelian extension of exponent $p$, i.e. the compositum
of cyclic extensions of order $p$. This means that
$\rho(\Gamma_{K_S})\subseteq V^\#\subseteq G^\#$.

Then Proposition \ref{central-prop} gives the result.
\end{pf}
\begin{prop}\label{socle-gl-n-k}
Let $K$, $n$, $p$, $E$, \O, \pG, $k$, $q$ and $S$ be as in Theorem \ref{theo}
and its proof. Let
$\rho_1$,~$\rho_2:\Gamma_K\longrightarrow\GL(n,E)$ be two representations
unramified outside $S$. Let $K'$ be the compositum of all extensions of $K$
unramified outside $S$ with degree $d$ such that:
\begin{itemize}
\item $d\mathop{|}\#\GL(n,k)$
\item $(d,p)=1$
\item $d<=\frac{q^n-1}{q-1}$.
\end{itemize}
Denote $\rho'_1$ and $\rho'_2$ the respective restrictions of $\rho_1$ and
$\rho_2$ to $\Gamma_{K'}$. Then $\rho'_1$ and $\rho'_2$ satisfy Condition
\ref{all-eigenvalues-equal} of Theorem \ref{theo}.
\end{prop}
\begin{pf}
Let $G$ be a subgroup of $\GL(n,k)$. We consider a flag
$V_0=\{0\}\subset V_1\subset\dots\subset V_\ell=k^n$ such that each $V_i$ is
stable under the action of $G$ and the action of $G$ on each quotient
$\overline{V_i}=\quot{V_i}{V_{i-1}}$ is irreducible. We will denote $G_i$ the
image of $G$ in $\Hom(\overline{V_i})$ and $d_i=\dim_k\overline{V_i}$. In a
basis adapted to the flag $(V_i)$, the matrices representing the action of $G$
on $k^n$ are blockwise upper-triangular and the $i$-th diagonal block of an
element $g\in G$ is equal to the projection of $g$ in $G_i$.

Then for any $i\in\{1,\dots,\ell\}$, $G_i$ is a finite group and a subgroup of
a general linear group. If $P$ is a $p$-Sylow subgroup of $G_i$, then the
elements of $P$ are the elements $g\in G_i$ such that, for certain basis
$(e_j)$ of $\overline{V_i}$, $\forall j\in\{1,\dots,d_i\}$,
$g(e_j)=e_j+\sum_{k<j}\lambda_{k,j}e_k$. In particular a $p$-Sylow of $G_i$
fixes at least one vector in $\overline{V_i}$. Let $e_i$ be such a vector,
$\{e_{i,j}\}_{1<=i<=n_i}$ its images under the action of $G_i$ (with
$e_{i,1}=e_i$) and
$H_{i,j}=\oper{Fix}_{G_i}(e_{i,j})=\{g\in G_i/g(e_{i,j})=e_{i,j}\}$. We see
that $H_{i,j}$ is a conjugate of $H_{i,1}=\oper{Fix}(e_i)$ and thus contains a
$p$-Sylow of $G_i$, in particular its index in $G_i$ is prime to $p$. Let
$H_i=\bigcap_jH_{i,j}$. Since $\overline{V_i}$ is irreducible under the action
of $G_i$, for any $v\in\overline{V_i}$ there exists $(\lambda_j)\in k^{d_i}$
such that $v=\sum_j\lambda_je_{i,j}$. In particular $\forall g\in H_i$,
$g(v)=v$ which means that $H_i=\{1\}$.

Let $G_{i,j}$ be the inverse image of each $H_{i,j}$ in $G$. Since $G_{i,j}$ is
a subgroup of $G$, we have
$[G:G_{i,j}]\mathop|\#G$. Since $G_i$ is a projection of $G$, we also have
$[G:G_{i,j}]=[G_i:H_{i,j}]<=q^{d_i}-1$ and
$([G:G_{i,j}],p)=([G_i:H_{i,j}],p)=1$. It is clear that
$G_{i,j}\cap Z(\GL(n,k))=\{1\}$ because all the elements of $H_{i,j}$ have
eigenvalues equal to $1$. Consider $G'=k^\times G$ then we also have
$[G':G_{i,j}]<=q^{d_i}-1$ (because $G_{i,j}$ is the fixator of $e_{i,j}$ also
in $G'$). This means that $[G':k^\times G_{i,j}]<=\frac{q^{d_i}-1}{q-1}$. This
in turns implies that if $Z_0=Z(\GL(n,k))\cap G$ and $G'_{i,j}=Z_0G_{i,j}$,
then $[G':k^\times G_{i,j}]=[G:G'_{i,j}]$ so that
$[G:G'_{i,j}]<=\frac{q^{d_i}-1}{q-1}$. The intersection of all the $G_{i,j}$
project trivially in each $G_i$, which means that its elements have eigenvalues
equal to $1$; we thus see that any $g\in\bigcap_{i,j}G'_{i,j}$ have all its
eigenvalues equal.

To finish the proof, take $G=(\rho_1\times\rho_2)(\Gamma_K)$ acting
on $k^n\times k^n$. Observe that $k^n\times\{0\}$ and $\{0\}\times k^n$ are
both stabilized by $G$ so that all $d_i<=n$. The inverse image in $\Gamma_K$ of
$G'_{i,j}$ defines an extension that have the properties listed in the
hypothesis of the proposition. This means that if $K_1$ is their compositum,
then the elements of $(\rho_1\times\rho_2)(\Gamma_{K_1})$ have all their
eigenvalues equal thus, since $K_1\subseteq K'$, $\rho'_1$ and $\rho'_2$
verify Condition \ref{all-eigenvalues-equal} of Theorem \ref{theo}.
\end{pf}
\section{Numerical application}\label{comparison}
% Galois: x^64 - 16*x^61 - 96*x^60 + 144*x^59 + 640*x^58 + 1424*x^57 + 1184*x^56 - 18960*x^55 - 41760*x^54 + 1376*x^53 + 197184*x^52 + 686112*x^51 + 503136*x^50 - 361488*x^49 - 32684*x^48 - 422688*x^47 + 3328944*x^46 + 194144*x^45 - 9106992*x^44 + 12742688*x^43 - 13880240*x^42 - 2172064*x^41 + 42205032*x^40 - 81439424*x^39 + 70223264*x^38 + 5170976*x^37 - 112924176*x^36 + 181443744*x^35 - 120283616*x^34 - 73923872*x^33 + 288559592*x^32 - 363513856*x^31 + 215744096*x^30 + 79679200*x^29 - 318677792*x^28 + 319483168*x^27 - 79843680*x^26 - 217273248*x^25 + 333944272*x^24 - 161711328*x^23 - 190908864*x^22 + 496539520*x^21 - 579760224*x^20 + 422942592*x^19 - 146636736*x^18 - 98472864*x^17 + 232483000*x^16 - 266632896*x^15 + 254039136*x^14 - 234357888*x^13 + 215933024*x^12 - 190302336*x^11 + 152557600*x^10 - 108211328*x^9 + 67231888*x^8 - 36439104*x^7 + 17140160*x^6 - 6942400*x^5 + 2395872*x^4 - 691136*x^3 + 159168*x^2 - 26240*x + 2308
%Gap: [64,34]
%Sous-groupes cycliques: [5, 137, 13, 1], [7, 257, 7, 1], [11, 73, 19, 1], [17, 337, 17, 1], [23, 257, 23, 1], [31, 1]
\subsection{Short version}
In \cite{GeemTop:non-self-dual}, the authors give an example of two non
self-dual representations of $\Gal(\quot{\overline\QM}\QM)$ (one should note
that the representation coming from the automorphic side is only conjectural)
and show that they have equal trace for all primes from $3$ to $67$. We can
apply our result to their example. In our terms, we have $n=3$, $p=2$ (so
that $m=2$), $K=\QM$, $E=\QM_2[i]$ and $S=\{2\ZM,\infty\}$. We denote by
$\rho_1$ and $\rho_2$ the representations they compare. There are no degree
$3$ and $7$ extensions of \QM that ramify only in $S$ so that, according to
Proposition \ref{socle-gl-n-k}, Condition \ref{all-eigenvalues-equal} of the
theorem is verified. We made a script in \textsc{gp/pari} to search for the
extensions described in the construction of the field $\QM_S$. We found that
the final compositum is a degree $64$ field, which we denote $\QM_{(2)}$. In
the paper \cite{GeemTop:non-self-dual}, it is shown that the characteristic
polynomial of the image of a Frobenius element $\Frob_\pG$ depends only on its
trace. As a consequence, all the eigenvalues of $\rho_i(\Frob_\pG)$ are
determined by $\Tr\rho_i(\Frob_\pG)$. The eigenvalues of
$\rho_i(\Frob_\pG^{\,k})=\rho_i(\Frob_\pG)^{\,k}$ are powers of the eigenvalues
of $\rho_i(\Frob_\pG)$, so that the characteristic polynomial of all the
$\rho_i(\Frob_\pG^{\,k})$ are determined by $\Tr\rho_i(\Frob_\pG)$. This means
that we can restrict the comparison to the traces of the images of the elements
generating maximal cyclic subgroups. Thanks to \textsc{gp/pari}, we found a
list of primes $p$ such that any element of the Galois group of $\QM_{(2)}$
over \QM is (conjugate to) the power of a Frobenius element above $p$. This
list is $\{5,7,11,17,23,31\}$. The prime $3$ is not included just because of
the method (and the particular polynomial defining $\QM_{(2)}$) used. Observe
that all of the primes have already been checked in the paper
\cite{GeemTop:non-self-dual}.

Professor Luis Dieulefait, from Universitat de Barcelona, made me observe that
on page 400 of the aforementioned article, the authors note that the
geometric representation is absolutely irreducible, which means in particular
that it is equal to its semi-simplification. The remark applies obviously also
to the conjectural automorphic representation.
\begin{cor}
The representation and the tentative representation compared in
\cite{GeemTop:non-self-dual} are isomorphic.
\end{cor}
Professor Dieulefait also observed that what is said about $P_5$ in the
aforementioned article is also true for $P_7=X^3-(1+4i)X^2+7(1+4i)X-7^3$ (the
field generated by one root of $P_7$ is of degree $6$ over \QM, contains only
fourth roots of unity and it is immediate to see that no rational multiple of
$i$ is a root of $P_7$). This means that all the members of the family of
$\ell$-adic representations are absolutely irreducible, hence semi-simple.
\enlargethispage{\baselineskip}
\subsection{Longer version}
The script used above to look for $K_S$ computes the sequence of fields
$(K_i)$. At each step, it computes linearly independent Kummer extensions of
$K_i$ and takes their compositum. Since the ramification is rather limited, we
tried to detect early (i.e. before computing the compositum) whether an
extension is not a sub-extension of $K_{i+1}$. For that purpose, we used the
fact that the residual extensions are cyclic, therefore we could not have
residual extensions larger than $p^m=4$. At each step the residual extension is
easily computed using class-field theory. We determined that the beginning of
the field sequence is as follows: $K_1$ is of degree $4$ over \QM, $K_2$ of
degree $32$
and $K_3$ of degree $64$. Since $\Gal(\quot{K_3}{K_2})$ is of order
$2$, $\Gal(\quot{K_S}{K_2})$ is a cyclic subgroup of
$\Gal(\quot{K_S}K)$\footnote{$\Gal(\quot{K_S}{K_3})$ is the Frattini subgroup
of $\Gal(\quot{K_S}{K_2})$ and can be of index $2$ if and only if
$\Gal(\quot{K_S}{K_2})$ is cyclic of order some power of $2$, here at most
$4$.}. Instead of looking for quadratic extensions of $K_3$, we checked that
cyclic extensions of order $4$ of $K_2$ all had a too large residual degree,
which proved that $K_3=K_4$ and thus $K_S=K_3$\footnote{As an additional proof,
we checked with \textsc{gap} that there is no order $128$ group admitting such
a chain of Frattini quotients.}. One equation for the extension is
\begin{multline*}
x^{64} - 16x^{61} - 96x^{60} + 144x^{59} + 640x^{58} + 1424x^{57} + 1184x^{56} - 18960x^{55} - 41760x^{54} +	\\
1376x^{53} + 197184x^{52} + 686112x^{51} + 503136x^{50} - 361488x^{49} - 32684x^{48} - 422688x^{47} +		\\
3328944x^{46} + 194144x^{45} - 9106992x^{44} + 12742688x^{43} - 13880240x^{42} - 2172064x^{41} +		\\
42205032x^{40} - 81439424x^{39} + 70223264x^{38} + 5170976x^{37} - 112924176x^{36} + 181443744x^{35} -		\\
120283616x^{34} - 73923872x^{33} + 288559592x^{32} - 363513856x^{31} + 215744096x^{30} + 79679200x^{29} -	\\
318677792x^{28} + 319483168x^{27} - 79843680x^{26} - 217273248x^{25} + 333944272x^{24} - 161711328x^{23} -\!	\\
190908864x^{22} + 496539520x^{21} - 579760224x^{20} + 422942592x^{19} - 146636736x^{18} - 98472864x^{17} +\!	\\
232483000x^{16} - 266632896x^{15} + 254039136x^{14} - 234357888x^{13} + 215933024x^{12} -			\\
190302336x^{11} + 152557600x^{10} - 108211328x^9 + 67231888x^8 - 36439104x^7 + 17140160x^6 -			\\
6942400x^5 + 2395872x^4 - 691136x^3 + 159168x^2 - 26240x + 2308.
\end{multline*}
Its Galois group is identified in Gap's small group library as $[64,34]$.
Up to conjugacy, this group has $6$ maximal cyclic subgroups. We list them
below using the following convention: if a cyclic subgroup is
$\{1,g,g^2,\dots,g^k\}$, we write it as $(1,p_1,p_2,\dots,p_k)$ where $p_i$ is
a prime number such that there is a Frobenius element above $p_i$ that is equal
to $g^i$. The list is:
$$
(1, 5, 137, 13); (1, 7, 257, 7); (1, 11, 73, 19); (1, 17, 337, 17); (1, 23, 257, 23); (1, 31)
$$
The center of the group is a two element subgroup generated by $\Frob(337)$.

Since we have $K_3=K_4$ and there are no extensions degree $3$, $5$, $7$, $9$
or $15$ of \QM ramifying only in $S$, the discussion above applies also to
$n=4$.

Hence, to test for the isomorphism of semi-simplification of representations of
$\Gal(\quot{\overline\QM}\QM)$ of dimension $3$ or $4$ over any finite
extension of $\QM_2$ having $\FM_2$ as residual field, ramifying
only at $2$ and $\infty$, it is sufficient to either test
\begin{itemize}
\item the traces at primes $\{5,7,11,13,17,19,23,31,73,137,257,337\}$;
\end{itemize}
or
\begin{itemize}
\item the characteristic polynomials at primes $\{5,7,11,17,23,31\}$.
\end{itemize}
\enlargethispage{\baselineskip}

\end{document}